\documentclass[12pt]{article} 
\usepackage[cp1251]{inputenc}
\usepackage[ukrainian]{babel}
\usepackage{latexsym,amsfonts,amssymb,amsmath,epsfig}
\usepackage{graphicx,graphics,hhline,cite}
\usepackage{euscript}

\begin{document}
\sloppy

 УДК 517.5

 \textbf{Анатолій Сердюк\footnote{Відповідальний за листування, e-mail: serdyuk@imath.kiev.ua}} (Інститут математики НАН України)

 \textbf{Сорич Віктор} (Кам'янець-Подільський національний університет  імені Івана Огієнка)

 \textbf{Ніна Сорич} (Кам'янець-Подільський національний університет  імені Івана Огієнка)

\vskip3mm

 \textbf{Найкраще наближення тригонометричними поліномами класів згорток, породжених деякими лінійними комбінаціями періодичних ядер}

\vskip3mm

Для довільних нетривіальних лінійних комбінацій скінченного числа ядер Пуассона \mbox{$P_{q_i,\beta}(t)=\sum\limits^\infty_{k=0}{q^k_i{\cos \left(kt-\frac{\beta\pi}{2}\right)\ }},\ \beta\in {\mathbb R},$ $q_i\in (0,1),$} встановлено виконання умови Надя $N^*_n$ для всіх номерів $n$, починаючи з деякого номера $n_0$. Також для будь-якого $n\in {\mathbb N}$ доведено існування лінійних комбінацій  $m\ (m\in\mathbb{N}\setminus\{1\})$ ядер   Бернуллі $D_{r_i}(t)={\left(-1\right)}^{\frac{r_i-1}{2}}\frac{{\sin k }t}{k^{r_i}}, r_i=2l_i-1, l_i\in {\mathbb N},$ де $r_i\ne r_j$ при $i\ne j$  , та лінійних  комбінацій $m$ спряжених ядер Пуассона $P_{q_i,1}(t)=\sum\limits^{\infty }_{k=1}q^k_i{\sin k}t,$ $ q_i\in (0,1),$ де $q_i\ne q_j$ при $i\ne j$, таких, що задовольняють умову Нікольського $A^*_n$ і при цьому не задовольняють умову Надя $N^*_n$. Як наслідок, в кожному з перелічених  випадків знайдено точні значення найкращих наближень в середньому таких лінійних комбінацій тригонометричними поліномами порядку не вищого за $n-1$ та обчислено точні значення найкращих наближень класів згорток , породжених зазначеними лінійними комбінаціями, в метриках просторів $C$ і $L$.

\vskip3mm
\textbf{1. Вступ}. Нехай $L_{\infty }$ -- простір вимірних $2\pi$-періодичних і істотно обмежених функцій $\varphi $ із нормою ${\left\|\varphi\right\|}_{L_\infty}={\left\|\varphi\right\|}_\infty=\mathop{\rm esssup}\limits_{t}\left|\varphi(t)\right|$ , C -- простір неперервних на всій дійсній осі $2\pi$-періодичних функцій $\varphi$ із нормою  ${\left\|\varphi\right\|}_C=\max\limits_{t}\left|\varphi(t)\right|$, $L=L_1-$ простір $2\pi$-періодичних сумовних на $\left(0,2\pi\right)$ функцій $\varphi$ із нормою${\left\|\varphi\right\|}_L={\left\|\varphi\right\|}_1=\int\limits^{2\pi}_0{\left|\varphi(t)\right|}dt$.

    Нехай $f\in L$. Величину
\begin{equation} \label{GrindEQ__1_} 
	E_n(f)_L=\inf_{t_{n-1}\in \tau_{2n-1}}{\left\|f-t_{n-1}\right\|}_1    ,                                          
\end{equation} 
в якій   inf  розглядають по множині $\tau_{2n-1}-$ усіх тригонометричних поліномів порядка $n-1,n\in \mathbb{N}$, називають найкращим наближенням в середньому функції $f$. Аналогічно для $f\in C$ величину
\begin{equation} \label{GrindEQ__2_} 
	E_n(f)_C=\inf_{t_{n-1}\in \tau_{2n-1}}{\left\|f-t_{n-1}\right\|}_C 
\end{equation} 
називають найкращим рівномірним наближенням функції  $f$.

    Нехай $K\in L$ і $\mathfrak{M}$ -- деяка підмножина функцій $\varphi$ із $L$. Через $\mathfrak{M}*K$ будемо позначати множину функцій $f$, які можна представити у вигляді згортки з ядром $K$, тобто у вигляді 
\begin{equation} \label{GrindEQ__3_} 
	f(x)=\frac{a_0}{2}+\frac{1}{\pi}\int\limits^{2\pi}_0{\varphi(x-t)K(t)dt=\frac{a_0}{2}}+\left(\varphi*K\right)(x),\ \varphi\in \mathfrak{M},a_0\in {\mathbb R}. 
\end{equation} 
Функції вигляду
\begin{equation} \label{GrindEQ__4_} 
	P_{q,\beta}(t)=\sum\limits^\infty_{k=1}{q^k{\cos \left(kt-\frac{\beta\pi}{2}\right)}}, q\in (0,1),\ \beta\in {\mathbb R}, 
\end{equation} 
називають ядрами Пуассона (див., наприклад, [1, c. 302], а функції вигляду
\begin{equation} \label{GrindEQ__5_} 
	D_r(t)=\sum\limits^\infty_{k=1}{\frac{1}{k^r}{\cos \left(kt-\frac{r\pi}{2}\right)}}, r\in {\mathbb N}, 
\end{equation} 
називають ядрами Бернуллі (див., наприклад, [1, c.137]).

    В рамках даної роботи в ролі $\mathfrak{M}$ будуть виступати множини       
\begin{equation} \label{GrindEQ__6_} 
	U^0_\infty=\left\{\varphi\in L_\infty:{\left\|\varphi\right\|}_\infty\le1,\int\limits^{2\pi}_0{\varphi(t)dt=0}\right\} 
\end{equation} 
та
\begin{equation} \label{GrindEQ__7_} 
	U^0_1=\left\{\varphi\in L_1:{\left\|\varphi\right\|}_1\le1,\int\limits^{2\pi}_0{\varphi(t)dt=0}\right\},                                          
\end{equation} 
а в якості твірних ядер $K$ розглядатимуться нетривіальні лінійні комбінації   \begin{equation} \label{GrindEQ__8_}
	K_m(t)=K_m(\overline{\alpha},t):=\sum\limits^m_{i=1}{\alpha_iK_i(t), \ \alpha_i}\in {\mathbb R},i=\overline{1,m},\sum\limits^m_{i=1}{\alpha^2_i>0,}
\end{equation}  
де $K_i(t)$ -- $2\pi$-періодичні функції, що є або ядрами Пуассона $P_{q,
\beta}(t)$, або ядрами Бернуллі $D_r(t),r=2l-1,l\in {\mathbb N}$.

    Для ядер $K_m(t)$ вигляду \eqref{GrindEQ__8_} вивчається задача про знаходження точних значень найкращих наближень в середньому, тобто величин
\begin{equation} \label{GrindEQ__9_} 
	E_n(K_m)_L=\inf_{t_{n-1}\in \tau_{2n-1}}{\left\|\sum\limits^m_{i=1}{\alpha_iK_i(\cdot)-t_{n-1}(\cdot)}\right\|}_1 
\end{equation} 
а також задача про точні значення величин
\begin{equation} \label{GrindEQ__10_} 
	E_n{\left(K_m*U^0_\infty\right)}_C=\sup_{f\in K_m*U^0_\infty}E_n(f)_C 
\end{equation} 
та
\begin{equation} \label{GrindEQ__11_} 
	E_n{\left(K_m*U^0_1\right)}_L=\sup_{f\in K_m*U^0_1}E_n(f)_L  ,                                                  
\end{equation} 
де $E_n(f)_L$ і $E_n(f)_C$  означені формулами \eqref{GrindEQ__1_} і \eqref{GrindEQ__2_} відповідно.

    Зазначимо, що у випадку, коли $m=1$ і $K(t)=K_m(t)=\alpha_1D_r(t),r\in {\mathbb N}{\rm ,}$ де  $D_r(t)-$ ядро Бернуллі вигляду \eqref{GrindEQ__5_}, класи $K*U^0_p=\alpha_1D_r*U^0_p,p=1,\infty ,$ з точністю до сталого множника $\alpha_1\neq0$ є відомими класами $W^r_p\ $ $2\pi$-періодичних функцій, що мають абсолютно неперервні похідні $f^{(k)}$до $r-1$  порядку включно і такі, що $f^{(r)}\in U^0_p$. Аналогічно, при $m=1$ і $K(t)=K_m(t)=\alpha_1P_{q,\beta}(t),q\in (0,1),$ де $P_{q,\beta}(t)-$ ядро Пуассона вигляду \eqref{GrindEQ__4_},  класи $K*U^0_p=\alpha_1P_{q,\beta}*U^0_p,p=1,\infty,$ з точністю до сталого множника $\alpha_1\neq0$ є відомими класами $C^q_{\beta,p} $ -- множин інтегралів Пуассона від функцій $\varphi$ з $U^0_p$ , тобто функцій $f$ вигляду
\begin{equation} \label{GrindEQ__12_} 
	f(x)=\frac{a_0}{2}+\frac{1}{\pi}\int\limits^\pi_{-\pi}{P_{q,\beta}}(x-t)\varphi(t)dt,\ \varphi\in U^0_p. 
\end{equation} 

Інтерес до твірних ядер $K$, зображуваних у вигляді лінійних комбінацій вигляду \eqref{GrindEQ__8_}, виникає в зв'язку з появою класів функцій, що породжуються диференціальними поліноміальними операторами зі сталими коефіцієнтами (див., наприклад, роботу М.Г. Крейна [2], де розглядались диференціальні оператори $P(D)=D^p+A_1D^{p-1}+...+A_p,D=\frac{d}{dt},$ роботу Н.І. Ахієзера [3], де були розглянуті оператори $L(D)=\prod^m_{k=1}{(D^2}+\alpha_k)$ та $M(D)=D\prod^m_{k=1}{(D^2+\alpha_k)}$, де  $\alpha_k\in {\mathbb R}$ та ін.).

    При  $m=1$ розв'язки задач \eqref{GrindEQ__9_}, \eqref{GrindEQ__10_} і \eqref{GrindEQ__11_} є відомими і у випадку, коли $K_1(t)=\alpha_1D_r(t),$ $ r\in {\mathbb N}, \alpha_1\neq0,$ випливають із робіт Ж. Фавара [4,5], 
 Н.І. Ахієзера та М.Г. Крейна [6], С.М.~Нікольського [7], а у випадку, коли $K_1(t)=\alpha_1 P_{q,\beta}(t),q\in (0,1), \alpha_1\neq0$ -- з робіт М.Г. Крейна [2], Б. Надя [8],
С.М. Нікольського [7] при $\beta\in {\mathbb Z}$ та із робіт А.В. Бушанського [9], В.Т. Шевалдіна [10] при $\beta\in {\mathbb R}\setminus {\mathbb Z}$.

    В усіх згаданих роботах обчислення точних значень величин вигляду \eqref{GrindEQ__9_}-\eqref{GrindEQ__11_} базувалось на встановленні того факту, що ядра Бернуллі $D_r(t),r\in {\mathbb N}$, та ядра Пуассона $ P_{q,\beta}(t),q\in (0,1),$ $\beta\in {\mathbb R},$ задовольняють так звану умову Нікольського $A^*_n$ і навіть більш жорстку, але і більш прозору, ніж $A^*_n$, умову Надя $N^*_n$.

 \textbf{Означення 1. }Кажуть, що сумовна $2\pi$-періодична функція $K(t)$, яка тотожно не дорівнює нулю, задовольняє умову$A^*_n, n\in {\mathbb N}, {\rm (}K\in A^*_n)$, якщо існують тригонометричний поліном $T^*_{n-1}\left(t\right)$ степеня $n-1$ і додатне число $\lambda\le \frac{\pi}{n}$ таке, що якщо  $\Delta(t)=K(t)-T^*_{n-1}(t), \varphi^*(t)={\rm sign}\,\Delta(t)$, то майже для всіх $t$ виконується рівність $\varphi^*(t+\lambda)=-\varphi^*(t)$. 

 В [7] було показано, що в якості $\lambda$ завжди можна взяти число $\lambda^*=\frac{\pi}{n^*}$, де $n^*$ - натуральне число $\ge n$. 

    \textbf{Означення 2.} Кажуть, що сумовна $2\pi$-періодична функція $K(t)$, яка тотожно не дорівнює нулю, задовольняє умову $N^*_n, n\in {\mathbb N}, (K\in N^*_n)$, якщо існують тригонометричний поліном $T^*_{n-1}\left(t\right)$ степеня $n-1$ і точка $\xi\in [0,\frac{\pi}{n})$ такі, що різниця $K(t)-T^*_{n-1}(t)$ змінює знак на $[0,2\pi)$ у точках $t_k=\xi+\frac{k\pi}{n},k=0,1,...,2n-1,$  і лише в них.  

    Із означень 1 і 2 безпосередньо випливає включення $N^*_n\subset A^*_n,n\in {\mathbb N}$. Водночас неважко привести приклади функцій $K$, поліномів $T^*_{n-1}\in \tau_{2n-1}$ і чисел $\xi, 0\le \xi<\frac{\pi}{n}$, для яких різниця $K(t)-T^*_{n-1}(t)$ почережно змінює знак більш ніж у $2n$ рівномірно розташованих точках на $[0,2\pi)$. В зв'язку з цією обставиною, пропонуємо наступне узагальнення умови $N^*_n$.

    \textbf{Означення 3.} Будемо казати , що сумовна $2\pi$-періодична функція $K(t)$, яка тотожно не дорівнює нулю, задовольняє умову$N^*_{n,p},n\in {\mathbb N}, p=0,1,... (K\in N^*_{n,p})$ якщо існують тригонометричний поліном $T^*_{n-1}\left(t\right)$ степеня $n-1$ і точка $\xi\in [0,\frac{\pi}{n+p})$ такі, що різниця $K(t)-T^*_{n-1}(t)$ змінює знак на $[0,2\pi)$ у точках $t_k=\xi+\frac{k\pi}{n+p},k=0,1,...,2n+2p-1,$  і лише в них.  

 Із означень 1-3 випливає, що при $p=0$ умова $N^*_{n,p}$ збігається з умовою Надя $N^*_n$ (тобто $N^*_{n,0}=N^*_n$ ) і, крім того, при всіх $p=0,1,...$ та $n\in {\mathbb N}$  мають місце включення $N^*_{n,p}\subset A^*_n$.

    Із результатів С.М. Нікольського [7] випливає, що для класів згорток $K*U_{\infty }$ і $K*U^0_1$, породжених ядрами $K(t)=K_m(t)$ вигляду \eqref{GrindEQ__8_}, які задовольняють умову $A^*_n$, виконуються наступні рівності:
\begin{equation} \label{GrindEQ__13_} 
	E_n(K*U^0_\infty)_C=\frac{1}{\pi}E_n(K)_L=\frac{1}{\pi}\int\limits^{2\pi}_0{\left|K(t)-T^*_{n-1}(t)\right|}dt, 
\end{equation} 
\begin{equation} \label{GrindEQ__14_} 
	E_n(K*U^0_1)_L=\frac{1}{\pi}E_n(K)_L=\frac{1}{\pi}\int\limits^{2\pi}_0{\left|K(t)-T^*_{n-1}(t)\right|}dt, 
\end{equation} 
в яких $T^*_{n-1}(t)-$ тригонометричний поліном, що фігурує в означенні 1.

    У випадку, коли ядро $K(t)$ задовольняє умову $N^*_{n,p}$, то крім \eqref{GrindEQ__13_}, \eqref{GrindEQ__14_} справджується рівність
\begin{equation} \label{GrindEQ__15_} 
	\frac{1}{\pi}\int\limits^{2\pi}_0{\left|K(t)-T^*_{n-1}(t)\right|}dt=\left|\frac{1}{\pi}\int\limits^{2\pi}_0{K(t)\, {\rm sign\,}{\sin \left((n+p)(t-\xi)\right)\ }dt}\right|, 
\end{equation} 
де $T^*_{n-1}(t)$ і $\xi$ -- тригонометричний поліном і точка, які фігурують в означенні 3.

 Зазначимо також, що завдяки роботам Б. Надя [8] та С.М. Нікольського [7] у випадку, коли твірні ядра $K$ задовольняють умову $N^*_n$ або $A^*_n$ можна вказати лінійний поліноміальний метод наближення $ $$L^*_{n-1}(f,)$ такий, що 
\[
\mathcal{E}{\left(K*U^0_\infty;L^*_{n-1}\right)}_C=\sup_{f\in K*U^0_\infty}{\left\|f(\cdot )-L^*_{n-1}(f;\cdot )\right\|}_C=E_n(K*U^0_\infty)_C,\] 
\[
\mathcal{E}{\left(K*U^0_1;L^*_{n-1}\right)}_L=\sup_{f\in K*U^0_1}{\left\|f(\cdot )-L^*_{n-1}(f;\cdot )\right\|}_L=E_n(K*U^0_1)_L.\] 
Останні рівності означають, що лінійний поліноміальний метод  $L^*_{n-1}$ є найкращим для класів згорток $K*U^0_{\infty }$ і $K*U^0_1$ серед усіх поліноміальних лінійних методів наближення в метриках просторів $C$ і $L$ відповідно.

 Крім того, (див., наприклад, [11-17]) виявилося, що підпростори тригонометричних поліномів є екстремальними підпросторами для колмогоровських, лінійних та бернштейнівських поперечників в просторах $C$ та $L$ класів згорток
 $K*U^0_{\infty }$ та $K*U^0_1$, породжених ядрами $K=D_r$, $K= P_{q,\beta}$ та деякими іншими.

    Зазначимо також, що інтерес до знаходження точних значень величин \eqref{GrindEQ__9_}-\eqref{GrindEQ__11_} при \mbox{$m\ge2$} мотивується тією обставиною, що до них зводиться знаходження розв'язків деяких інших екстремальних задач теорії апроксимації, таких, зокрема, як задача Колмогорова-Нікольського по наближенню класів згорток полігармонічними ядрами Пуассона (див., наприклад, роботи [18-19] і наявну там бібліографію).

    В даній роботі для довільних нетривіальних лінійних комбінацій \eqref{GrindEQ__8_}, що породжуються ядрами Пуассона $P_{q_i,\beta}(t)$ вигляду \eqref{GrindEQ__4_}, $\beta\in {\mathbb R}, q_i\in (0,1), q_i\neq q_j$  при $i\ne j$, встановлено виконання умови Надя $N^*_n$ для всіх номерів $n$, починаючи з деякого номера $n_0$ в, як наслідок, для зазначених $n$ обчислено точні значення величин \eqref{GrindEQ__9_}, \eqref{GrindEQ__10_} і \eqref{GrindEQ__11_}.

    Крім того, для будь-якого $n\in {\mathbb N}$ в роботі доведено існування лінійних комбінацій $m\left(m\in {\mathbb N}{\rm \backslash }\left\{1\right\}\right)$ ядер Бернуллі $D_{r_i}(t),r_i=2l_i-1,l_i\in {\mathbb N}{\rm ,}$ та лінійних комбінацій $m\left(m\in {\mathbb N}{\rm \backslash }\left\{1\right\}\right)$ ядер Пуассона $P_{q_i,1}(t)$,$q_i\in (0,1),q_i\ne q_j$при $i\ne j$, таких, що задовольняють умову $N^*_{n,m-1}$. Як наслідок, для зазначених лінійних комбінацій знайдено точні значення величин \eqref{GrindEQ__9_}, \eqref{GrindEQ__10_} і \eqref{GrindEQ__11_}.

\vskip3mm
 \textbf{2. Основні результати. }

    Розглянемо випадок, коли ядро $K$, що породжує класи згорток, має вигляд 
\begin{equation} \label{GrindEQ__16_} 
	K(t)=K_m(t)=P^{\overline{\alpha},\overline{q}}_{m,\beta}(t)=\sum\limits^m_{i=1}{\alpha_i\sum\limits^\infty_{k=1}{q^k_i{\cos \left(kt-\frac{\beta\pi}{2}\right)}}},                    
\end{equation} 
\begin{equation} \label{GrindEQ__17_} 
	m\in {\mathbb N}, \alpha_i\in {\mathbb R}, \alpha_1\neq0,0<q_m<q_{m-1}<...<q_2<q_1<1.                
\end{equation} 
Зрозуміло, що довільна нетривіальна лінійна комбінація скінченного числа різних ядер Пуассона вигляду \eqref{GrindEQ__4_} може бути подана у вигляді \eqref{GrindEQ__16_} при деяких 
 наборах $\overline{\alpha}={\left\{\alpha_i\right\}}^m_{i=1}$ та $\overline{q}={\left\{q_i\right\}}^m_{i=1}$, для котрих виконуються умови \eqref{GrindEQ__17_}. 

    \textbf{Теорема 1}. Для будь-якого ядра $P^{\overline{\alpha},\overline{q}}_{m,\beta}(t)$ вигляду \eqref{GrindEQ__16_} за виконання умов \eqref{GrindEQ__17_} знайдеться номер $n_0$ такий, що для всіх $n\ge n_0, n\in {\mathbb N}$, справджується включення  $P^{\overline{\alpha},\overline{q}}_{m,\beta}\in N^*_n$ і, як наслідок, мають місце рівності
\[
E_{n} \left(P_{m,\beta }^{\overline{\alpha },\overline{q}} *U_{\infty }^{0} \right)_{C} =E_n{\left(P^{\overline{\alpha},\overline{q}}_{m,\beta}*U^0_1\right)}_L=\frac{1}{\pi}E_n{\left(P^{\overline{\alpha},\overline{q}}_{m,\beta}\right)}_L
={\left\|P^{\overline{\alpha},\overline{q}}_{m,\beta}*{\rm sign\ }{\sin n}(\cdot )\right\|}_C=
\] 
\begin{equation} \label{GrindEQ__18_} 
	=\frac{4}{\pi}\left|\sum\limits^m_{i=1}{\alpha_i\sum\limits^\infty_{k=0}{\frac{q^{(2k+1)n}_i}{2k+1}{\sin \left((2k+1)\theta_n\pi-\frac{\beta\pi}{2}\right)\ }}}\right|, 
\end{equation} 
де  $\theta_n$ -- єдиний на $[0,1)$ корінь рівняння
\begin{equation} \label{GrindEQ__19_} 
	\sum\limits^m_{i=1}{\alpha_i}\sum\limits^\infty_{k=0}{q^{(2k+1)n}_i}{\cos \left((2k+1)\theta_n\pi-\frac{\beta\pi}{2}\right)\ }=0.                                    
\end{equation} 
\textbf{Доведення}. Запишемо ядро $K_m(t)$ в \eqref{GrindEQ__16_} у вигляді
\begin{equation} \label{GrindEQ__20_} 
	K(t)=K_m(t)=\Psi_\beta(t)=\sum\limits^\infty_{k=1}{\psi(k){\cos \left(kt-\frac{\beta\pi}{2}\right)}},                              
\end{equation} 
де 
\begin{equation} \label{GrindEQ__21_} 
	\psi(k):=\sum\limits^m_{i=1}{\alpha_iq^k_i}  .                                                              
\end{equation} 
Покажемо, що послідовність  $\psi(k)$ вигляду \eqref{GrindEQ__21_} за умов \eqref{GrindEQ__17_} задовольняє умову Даламбера $\mathcal{D}_q, q\in (0,1)$ при $q=q_1$, тобто
\begin{equation}\label{GrindEQ__22_} 
	\lim_{k\to \infty }\frac{\psi(k+1)}{\psi(k)}=q_1,
\end{equation} 
 цей факт будемо коротко записувати $\psi\in \mathcal{D}_{q_1}$.  Дійсно, згідно \eqref{GrindEQ__21_} і \eqref{GrindEQ__17_}
$$
\frac{\psi(k+1)}{\psi(k)}=\frac{\sum\limits^m_{i=1}{\alpha_iq^{k+1}_i}}{\sum\limits^m_{i=1}{\alpha_iq^k_i}}=
\frac{\alpha_1q^{k+1}_1+\sum\limits^m_{i=2}{\alpha_iq^{k+1}_i}}{\alpha_1q^k_1+\sum\limits^m_{i=2}{\alpha_iq^k_i}}=
\frac{q_1+\sum\limits^m_{i=2}{\frac{\alpha_i}{\alpha_1}{\left(\frac{q_i}{q_1}\right)}^kq_i}}{1+\sum\limits^m_{i=2}{\frac{\alpha_i}{\alpha_1}{\left(\frac{q_i}{q_1}\right)}^k}}=q_1+\delta_k,$$  
де
\begin{equation} \label{GrindEQ__23_} 
	\delta_k:=\frac{\psi(k+1)}{\psi(k)}-q_1=\frac{\sum\limits^m_{i=2}{\frac{\alpha_i}{\alpha_1}{\left(\frac{q_i}{q_1}\right)}^k(q_i-q_1)}}{1+\sum\limits^m_{i=2}{\frac{\alpha_i}{\alpha_1}{\left(\frac{q_i}{q_1}\right)}^k}}  .                                 
\end{equation} 
Оскільки в силу \eqref{GrindEQ__17_} $0<\frac{q_i}{q_1}<1, i=2,3,...,m,$ то як випливає з \eqref{GrindEQ__23_}
\begin{equation} \label{GrindEQ__24_} 
	\lim_{k\rightarrow\infty}\delta_k=0  .                                                             
\end{equation} 
Рівність \eqref{GrindEQ__24_} доводить \eqref{GrindEQ__22_}.

    Згідно з теоремою 5 роботи [20] (див. також [21, 22]) якщо $\psi\in \mathcal{D}_q,q\in (0,1)$, то знайдеться залежний від $\psi$ номер $n_0$ такий, що для будь-якого натурального числа $n$, більшого або рівного за $n_0$, має місце включення $K=\Psi_\beta\in N^*_n$ і виконуються рівності \eqref{GrindEQ__13_}, \eqref{GrindEQ__14_}, а також рівність \eqref{GrindEQ__15_} при $p=0$. При цьому  $\xi=\theta_n\pi/n$ , де $\theta_n$ -- єдиний на $[0,1)$ корінь рівняння 
\begin{equation} \label{GrindEQ__25_} 
	\sum\limits^\infty_{k=0}{\psi\left((2k+1)n\right)}{\cos \left((2k+1)\theta_n\pi-\frac{\beta\pi}{2}\right)\ }=0,                             
\end{equation} 
a поліном $T^*_{n-1}$ однозначно визначається інтерполяційними умовами
\[
T^*_{n-1}\left(\frac{\theta_n\pi+k\pi}{n}\right)=\Psi_\beta\left(\frac{\theta_n\pi+k\pi}{n}\right),k=0,1,...,2n-2.
\] 
При $\psi(k)$, що мають вигляд \eqref{GrindEQ__21_} із \eqref{GrindEQ__13_}-\eqref{GrindEQ__15_} випливають рівності \eqref{GrindEQ__19_}. Теорему доведено.

    \textbf{Зауваження 1}. Як випливає з доведення теореми 3 роботи [22], номер $n_0$, який фігурує в теоремі 1, може бути пред'явлений у явному вигляді, як найменше натуральне число $n$, для якого 
\begin{equation} \label{GrindEQ__26_} 
	(1-q_1)^2\ge\frac{5+3q^2_1}{1-q^2_1}\cdot \frac{(1+q^2_1)^{2n}}{2^{2n}\sqrt{1-(1+q^2_1)^{2n}/2^{2n-1}}}+\varepsilon_n(2+\varepsilon_n),n=n_0,n_0+1,..., 
\end{equation} 
де $\varepsilon_n:=\sup\limits_{k\ge n}\left|\delta_k\right|$ , а $\delta_k$ означені в \eqref{GrindEQ__23_}. Існування такого номера $n_0$ випливає з факту монотонного спадання до нуля правої частини нерівності \eqref{GrindEQ__26_} по $n$ при $n\rightarrow\infty$.

    При $\beta\in {\mathbb Z}$ рівності \eqref{GrindEQ__18_} можуть бути конкретизовані. Якщо $\beta=2\nu, \nu\in {\mathbb Z}$ , то неважко збагнути, що для довільних наборів ${\left\{\alpha_i\right\}}^m_{i=1}$ і ${\left\{q_i\right\}}^m_{i=1}$ , підпорядкованих умовам \eqref{GrindEQ__17_}, єдиним на $[0,1)$ коренем рівняння  \eqref{GrindEQ__19_} є значення $\theta_n=\frac{1}{2}$  і в цьому випадку
\[
\sum\limits^m_{i=1}{\alpha_i\sum\limits^\infty_{k=0}{\frac{q^{(2k+1)n}_i}{2k+1}}}{\sin \left((2k+1)\theta_n\pi-\frac{\beta\pi}{2}\right)\ }=(-1)^\nu\sum\limits^m_{i=1}{\alpha_i\sum\limits^\infty_{k=0}{(-1)^k\frac{q^{(2k+1)n}_i}{2k+1}}}=
\]
\begin{equation} \label{GrindEQ__27_}
	=(-1)^\nu\sum\limits^m_{i=1}{\alpha_i\arctg q^n_i}  .                                                                          
\end{equation} 
Аналогічно, якщо $\beta=2\nu-1, \nu\in {\mathbb Z}$, то єдиним  на $[0,1)$ коренем рівняння \eqref{GrindEQ__19_} є значення $\theta_n=0$  і тоді 
\[
\sum\limits^m_{i=1}{\alpha_i\sum\limits^\infty_{k=0}{\frac{q^{(2k+1)n}_i}{2k+1}}}{\sin \left((2k+1)\theta_n\pi-\frac{\beta\pi}{2}\right)\ }=(-1)^\nu\sum\limits^m_{i=1}{\alpha_i\sum\limits^\infty_{k=0}{\frac{q^{(2k+1)n}_i}{2k+1}}}=\] 
\begin{equation} \label{GrindEQ__28_} 
	=(-1)^\nu\sum\limits^m_{i=1}{\frac{\alpha_i}{2}}{\ln \frac{1+q^n_i}{1-q^n_i}\ }.                                                                  
\end{equation} 
На підставі \eqref{GrindEQ__27_} і \eqref{GrindEQ__28_} при $\beta\in {\mathbb Z}$ із теореми 1 одержуємо таке твердження.

    \textbf{Теорема 2}. Якщо $\beta\in {\mathbb Z}$, то за виконання умов \eqref{GrindEQ__17_} ядра $P^{\overline{\alpha},\overline{q}}_{m,\beta}$  вигляду \eqref{GrindEQ__16_} задовольняють умову Надя $N^*_n$ при всіх номерах $n\in {\mathbb N}$, починаючи з деякого номера $n_0$, і для всіх таких $n$ мають місце рівності
\[E_{n} \left(P_{m,\beta}^{\overline{\alpha },\overline{q}} *U_{\infty }^{0} \right)_{C} =E_n{\left(P^{\overline{\alpha},\overline{q}}_{m,\beta}*U^0_1\right)}_L=\frac{1}{\pi}E_n{\left(P^{\overline{\alpha},\overline{q}}_{m,\beta}\right)}_L=\] 
\begin{equation} \label{GrindEQ__29_} 
	={\left\|P^{\overline{\alpha},\overline{q}}_{m,\beta}*{\rm sign\ }{\sin n}(\cdot )\right\|}_C=\left\{ \begin{array}{l}
		\frac{4}{\pi}\left|\sum\limits^m_{i=1}{\alpha_i \arctg q^n_i}\right|,\ \beta=2\nu, \nu\in {\mathbb Z}, \\ 
		\frac{2}{\pi}\left|\sum\limits^m_{i=1}{\alpha_i{\ln \frac{1+q^n_i}{1-q^n_i}\ }}\right|,\ \beta=2\nu-1, \nu\in {\mathbb Z}. \end{array}
	\right. 
\end{equation} 
Зазначимо, що номер $n_0$, який фігурує у теоремах 1 і 2 залежить від $m$, наборів  ${\left\{\alpha_i\right\}}^m_{i=1}$ та ${\left\{q_i\right\}}^m_{i=1}$ і не залежить від параметра $\beta$. Далі ми покажемо, що при кожному фіксованому $n\in {\mathbb N}$ існують лінійні комбінації $K_m$ різних ядер Бернуллі ${\left\{D_{r_i}\right\}}^m_{i=1},r_i=2l_i-1,l_i\in {\mathbb N}$ , а також лінійні комбінації $K_m$ різних спряжених ядер Пуассона ${\left\{P_{q_i,1}\right\}}^m_{i=1},q_i\in \left(0,1\right)$, такі, що задовольняють умову $N^*_{n,m-1}$ при $m>1$ (тобто $K_m\notin N^*_n$ і, водночас, $K_m\in N^*_{n,m-1}\subset A^*_n$).

    Мають місце наступні твердження.

    \textbf{Теорема 3}. Нехай $n\in {\mathbb N}$, $m\in {\mathbb N}{\rm \backslash }\left\{1\right\}$, і задано довільний  впорядкований набір $\overline{r}=(r_1,r_2,...,r_m)$ $m$ різних непарних чисел ($r_i=2l_i-1,l_i\in {\mathbb N},i=\overline{1,m}$).  Тоді існує впорядкований набір ${\overline{\alpha}}^*=(\alpha^*_1,\alpha^*_2,...,\alpha^*_m)$ дійсних чисел $\alpha^*_i, \sum\limits^m_{i=1}{(\alpha^*_i)}^2>0,$ та існує тригонометричний поліном $T^*_{n-1}$ степеня не вищого за $n-1$ , який інтерполює лінійну комбінацію 
\begin{equation} \label{GrindEQ__30_} 
	K_m(t)=D^{{\overline{\alpha}}^*,\overline{r}}_{m,1}(t)=\sum\limits^m_{i=1}{\alpha^*_i}D_{r_i}(t)=\sum\limits^m_{i=1}{\alpha^*_i}(-1)^{\frac{r_i-1}{2}}\sum\limits^\infty_{k=1}{\frac{{\sin k}t}{k^{r_i}}} 
\end{equation} 
на інтервалі $\left(0,2\pi\right)$ лише в точках $t_k=\frac{k\pi}{n+m-1}, k=1,...,2n+2m-3,$ , в яких також різниця $K_m(t)-T^*_{n-1}(t)$ почережно змінює знак, тобто:
\begin{equation} \label{GrindEQ__31_} 
	{\rm sign}\left(K_m(t)-T^*_{n-1}(t)\right)=\delta\,{\rm sign}\sin (n+m-1)t,                                    
\end{equation} 
де  $\delta=\pm 1$.

    \textbf{Теорема 4.} Нехай $n\in {\mathbb N}$, $m\in {\mathbb N}{\rm \backslash }\left\{1\right\}$, і задано довільний вектор $\overline{q}=(q_1,q_2,...,q_m)$ з різних чисел $q_i\in (0,1),i=\overline{1,m}$ . Тоді знайдеться впорядкований набір ${\overline{\alpha}}^*=(\alpha^*_1,\alpha^*_2,...,\alpha^*_m)$, дійсних чисел $\alpha^*_i,\sum\limits^m_{i=1}{(\alpha^*_i})^2>0,$  та тригонометричний многочлен $T^*_{n-1}$ степеня не вищого за $n-1$, який інтерполює лінійну комбінацію 
\begin{equation} \label{GrindEQ__32_} 
	K_m(t)=K_m({\overline{\alpha}}^*,t)=P^{{\overline{\alpha}}^*,\overline{q}}_{m,1}(t)=\sum\limits^m_{i=1}{\alpha^*_iP_{q_i,1}}(t)=\sum\limits^m_{i=1}{\alpha^*_i}\sum\limits^\infty_{k=1}{q^k_i{\sin k\ }t} 
\end{equation} 
на періоді $[0,2\pi)$лише в точках $t_k=\frac{k\pi}{n+m-1},k=0,1,...,2n+2m-3$, в яких до того ж різниця $K_m(t)-T^*_{n-1}(t)$ почережно змінює знак, тобто виконується рівність \eqref{GrindEQ__31_}.

    Фактично теореми 3 і 4 стверджують, що при $n\in {\mathbb N}$, $m\in {\mathbb N}\setminus\left\{1\right\}$ лінійні комбінації $D^{{\overline{\alpha}}^*,\overline{r}}_{m,1}(t)$ та $P^{{\overline{\alpha}}^*,\overline{q}}_{m,1}(t)$ задовольняють умову $N^*_{n,m-1}$ при $\xi=\frac{\pi}{2(n+m-1)}$, а, отже, для них не виконується умова Надя $N^*_n$, але виконується умова Нікольського $A^*_n$. Тому в силу рівностей \eqref{GrindEQ__13_} і \eqref{GrindEQ__14_}, а також формули \eqref{GrindEQ__15_} при $p=m-1$ із теорем 3 і 4 одержуємо наступні твердження.

    \textbf{Теорема 5}. Нехай $n\in {\mathbb N}$, $m\in {\mathbb N}\setminus\left\{1\right\}$ і задано довільний  впорядкований набір $\overline{r}=(r_1,r_2,...,r_m)$ різних непарних чисел ($r_i=2l_i-1,l_i\in {\mathbb N}{\rm ,}i=\overline{1,m}$). Тоді існує впорядкований набір ${\overline{\alpha}}^*=(\alpha^*_1,\alpha^*_2,...,\alpha^*_m)$ дійсних чисел $\alpha^*_i, \sum\limits^m_{i=1}{(\alpha^*_i})^2>0,$ такий, що для лінійної комбінації $K_m(t)=D^{{\overline{\alpha}}^*,\overline{r}}_{m,1}(t)$ вигляду \eqref{GrindEQ__30_} виконуються рівності
\[
E_n{\left(D^{{\overline{\alpha}}^*,\overline{r}}_{m,1}*U^0_\infty\right)}_C=E_n{\left(D^{{\overline{\alpha}}^*,\overline{r}}_{m,1}*U^0_1\right)}_L=\frac{1}{\pi}E_n{\left(D^{{\overline{\alpha}}^*,\overline{r}}_{m,1}\right)}_L=\] 
\begin{equation} \label{GrindEQ__33_} 
	={\left\|D^{{\overline{\alpha}}^*,\overline{r}}_{m,1}*{\rm sign\ }{\sin (}n+m-1)(\cdot )\right\|}_C=\left|\sum\limits^m_{i=1}{\alpha^*_i\frac{(-1)^{\frac {r_i-1}2}K_{r_i}}{(n+m-1)^{r_i}}}\right|,                             
\end{equation} 
де $K_{r_i}-$ константи Фавара
\[K_{r_i}:=\frac{4}{\pi}\sum\limits^\infty_{k=0}{\frac{(-1)^{k(r_i-1)}}{(2k+1)^{r_i+1}}=\frac{4}{\pi}}\sum\limits^\infty_{k=0}{\frac{1}{(2k+1)^{r_i+1}}}.\] 
\textbf{Теорема 6.} Нехай $n\in {\mathbb N}$, $m\in {\mathbb N}{\rm \backslash }\left\{1\right\}$ і задано довільний вектор $\overline{q}=(q_1,q_2,...,q_m)$ з різних чисел $q_i\in (0,1),i=\overline{1,m}$ . Тоді існує впорядкований набір ${\overline{\alpha}}^*=(\alpha^*_1,\alpha^*_2,...,\alpha^*_m)$ дійсних чисел  $\alpha^*_i$ такий, що для лінійної комбінації $K_m(t)=P^{{\overline{\alpha}}^*,\overline{q}}_{m,1}(t)$ вигляду \eqref{GrindEQ__32_} виконуються рівності          
\[
E_n{\left(P^{{\overline{\alpha}}^*,\overline{q}}_{m,1}*U^0_\infty\right)}_C=E_n{\left(P^{{\overline{\alpha}}^*,\overline{q}}_{m,1}*U^0_1\right)}_L=\frac{1}{\pi}E_n{\left(P^{{\overline{\alpha}}^*,\overline{q}}_{m,1}\right)}_L=\left\|P^{{\overline{\alpha_i}}^*,\overline{q}}_{m,1}*{\rm sign\ }{\sin (}n+m-1)(\cdot )\right\|=
\]
\begin{equation} \label{GrindEQ__34_} 
	=\frac{4}{\pi}\left|\sum\limits^m_{i=1}{\alpha^*_i\sum\limits^\infty_{k=0}{\frac{q^{(2k+1)(n+m-1)}_i}{2k+1}}}\right|=\frac{2}{\pi}\left|\sum\limits^m_{i=1}{\alpha^*_i{\ln \frac{1+q^{n+m-1}_i}{1-q^{n+m-1}_i}\ }}\right|.                                     
\end{equation} 

\textbf{3. Доведення теорем 3 і 4.}

    Для доведення теорем 3 та 4 знадобиться наступне допоміжне твердження, яке, можливо, не позбавлене і самостійного значення.

    \textbf{Лема 1}. Нехай $n\in {\mathbb N}$ і  $K_i(t),i=\overline{1,m},$ -- набір непарних і неперервних на $\left(0,2\pi\right)$  $2\pi$-періодичних функцій. Тоді для будь-якого набору ${\left\{t_k\right\}}^{n+l-1}_{k=1}$ різних точок  $t_k\in \left(0,\pi\right),$ $ l\in \left\{1,2,...,m-1\right\}$, існує нетривіальна лінійна комбінація $K_m\left({\overline{\alpha}}^*,t\right)=\sum\limits^m_{i=1}{\alpha^*_iK_i}(t)$ та непарний тригонометричний многочлен $T^*_{n-1}(t)$, який інтерполює $K_m\left({\overline{\alpha}}^*,t\right)$ в точках $t_k\in (0;\pi),$ $k=\overline{1;n+l-1}.$

    \textbf{Доведення леми 1}. Твердження леми 1 є очевидним, якщо серед нетривіальних лінійних комбінацій функцій ${\left\{K_i\right\}}^m_{i=1}$ існує така, що перетворюється у тригонометричний поліном порядку $\le n-1$. В цьому випадку згадані лінійна комбінація і поліном є шуканими і до того ж $K_m\left({\overline{\alpha}}^*,t\right)-T^*_{n-1}(t)\equiv 0$. Цей випадок можна вважати тривіальним. Зрозуміло, що тривіальний випадок виникає, зокрема, коли функції з набору $K_i(t), i=\overline{1,m},$ є лінійно залежними. Припустимо, що будь-яка нетривіальна лінійна комбінація $\sum\limits^m_{i=1}{\alpha_iK_i}(t)$ не є тригонометричним поліномом порядку $\le n-1$. Оскільки будь-яку непарну і неперервну на $\left(0,2\pi\right)$ функцію можна проінтерполювати непарним  тригонометричним поліномом порядку $n-1$ за системою $n-1$ довільних наперед заданих різних точок (див., наприклад, [13, c.62]), то для довільного $i=\overline{1,m}$ через $T^{(i)}_{n-1}(t)$ позначимо тригонометричний многочлен порядку $n-1$, який інтерполює функцію $K_i(t)$ в точках  $t_k\in (0;\pi), k=\overline{1,n-1}$. Нехай  $c_{i,k}=K_i(t_k)-T^{(i)}_{n-1}(t_k)$ , тоді  при $k=\overline{1,n-1}\ c_{i,k}=0, i=\overline{1,m}.$ Розглянемо $l$-вимірні вектори $c^{(1)}=(c_{1,n};c_{1,n+1};...;c_{1,n+l-1}),$  $c^{(2)}=(c_{2,n};c_{2,n+1};...;c_{2,n+l-1}),...,$  $c^{(m)}=(c_{m,n};c_{m,n+1};...;c_{m,n+l-1})$. Оскільки $l<m$ , то дана система векторів  є лінійно залежною, тому існує їх нетривіальна лінійна комбінація, що рівна нуль-вектору $\theta$:    $\sum\limits^m_{i=1}{\alpha^*_ic^{(i)}}=\theta$ , тобто виконується наступна система рівностей:
\begin{equation} \label{GrindEQ__35_} 
	\left\{ \begin{array}{l}
		\sum\limits^m_{i=1}{\alpha^*_ic_{i,n}}=0, \\ 
		\sum\limits^m_{i=1}{\alpha^*_i}c_{i,n+1}=0, \\ 
		.......................... \\ 
		\sum\limits^m_{i=1}{\alpha^*_i}c_{i,n+l-1}=0, \end{array}
	\right. 
\end{equation} 
де  $\sum\limits^m_{i=1}{(\alpha^*_i})^2>0$ .                            

       Покладемо  $\alpha^*=(\alpha^*_1,\alpha^*_2,...,\alpha^*_m)$, $T^*_{n-1}(t)=\sum\limits^m_{i=1}{\alpha^*_iT^{(i)}_{n-1}}(t)$. Покажемо, що вектор $\alpha^*$ і многочлен $T^*_{n-1}(t)$ шукані. 

    Дійсно, при $k=\overline{1,n-1}$ маємо    $K_m(\overline{\alpha}^*,t_k)-T^*_{n-1}(t_k)=\sum\limits^m_{i=1}{\alpha^*_iK_i}(t_k)-\sum\limits^m_{i=1}{\alpha^*_iT^{(i)}_{n-1}}(t_k)=\sum\limits^m_{i=1}{\alpha^*_i}\left(K_i(t_k)-T^{(i)}_{n-1}(t_k)\right)=0,$ 
 а при $k=\overline{n,n+l-1}$ в силу  рівностей системи \eqref{GrindEQ__35_}  маємо $K_m(\overline{\alpha}^*,t_k)-T^*_{n-1}(t_k)=\sum\limits^m_{i=1}{\alpha^*_i}\left(K_i(t_k)-T^{(i)}_{n-1}(t_k)\right)=\sum\limits^m_{i=1}\alpha^*_ic_{i,k}=0$.   Отже,  в усіх точках $t_k\in (0,\pi), k=\overline{1,n+l-1},$ функція $K_m({\overline{\alpha}}^*,t)$   та многочлен $T^*_{n-1}(t)$ набувають однакових значень. Лему 1 доведено.

       При   $l=m-1$ із леми 1 випливає таке твердження.

    \textbf{Наслідок 1}. Нехай $n,m\in {\mathbb N}$ і  $K_i(t),i=\overline{1,m},$ -- набір непарних і неперервних на $\left(0,2\pi\right)$  $2\pi$-періодичних функцій. Тоді  існує нетривіальна лінійна комбінація $K_m\left({\overline{\alpha}}^*,t\right)=\sum\limits^m_{i=1}{\alpha^*_iK_i}(t)$ та непарний тригонометричний многочлен $T^*_{n-1}(t)$, який інтерполює $K_m\left({\overline{\alpha}}^*,t\right)$ на $\left(0,2\pi\right)$ в \mbox{$2(n+m)-3$} рівномірно розподілених точках  $t_k$ вигляду 
\begin{equation} \label{GrindEQ__36_} 
	t_k=\frac{k\pi}{n+m-1}, k=\overline{1,2(n+m)-3}. 
\end{equation} 
\textbf{Доведення наслідку 1}. При $n,m\in {\mathbb N}$ і  $l=m-1$ для рівномірно розподілених точок $t_k=\frac{k\pi}{n+m-1}$ згідно з лемою 1 знайдеться нетривіальна лінійна комбінація  $\sum\limits^m_{i=1}{\alpha^*_iK_i(t)}$ і непарний тригонометричний многочлен $T^*_{n-1}(t)$ такі, що при всіх $k=\overline{1,n+m-2}$
\begin{equation} \label{GrindEQ__37_} 
	K_m\left({\overline{\alpha}}^*,\frac{k\pi}{n+m-1}\right)=T^*_{n-1}\left(\frac{k\pi}{n+m-1}\right).                             \end{equation} 

Із неперності і $2\pi$-періодичності функцій $K_m({\overline{\alpha}}^*,t)$  і  $T^*_{n-1}(t)$ випливає, що рівності \eqref{GrindEQ__37_} мають місце також і для $k=\overline{n+m,2(n+m)-3}$. Із непарності та неперервності функції $K_m\left({\overline{\alpha}}^*,t\right)-T^*_{n-1}(t)$ випливає також, що $K_m\left({\overline{\alpha}}^*,\pi\right)=T^*_{n-1}(\pi)$, тобто виконання \eqref{GrindEQ__37_} і при $k=n+m-1$. Наслідок 1 доведено.

    \textbf{Доведення теореми 3.} Оскільки непарні ядра Бернуллі   
    \begin{equation} \label{GrindEQ__38_} D_{r_i}(t)=\sum\limits^\infty_{k=1}{\frac{{\cos \left(kt-\frac{r_i\pi}{2}\right)\ }}{k^{r_i}}}={\left(-1\right)}^{\frac{r_i-1}{2}}\sum\limits^\infty_{k=1}{\frac{{\sin k }t}{k^{r_i}}}, r_i=2l_i-1, l_i\in {\mathbb N}, i=\overline{1,m},
    \end{equation} 
є неперервними на $\left(0,2\pi\right)$ функціями, то відповідно до наслідку 1 існує нетривіальна лінійна комбінація $K_m({\overline{\alpha}}^*,t)$ вигляду \eqref{GrindEQ__30_} і непарний тригонометричним поліномом $T^*_{n-1}(t)$ порядку не вищого за $n-1$, який інтерполює $K_m({\overline{\alpha}}^*,t)$ на $\left(0,2\pi\right)$ в $2(n+m)-3$ точках вигляду \eqref{GrindEQ__36_}.

    Далі ми покажемо, що жоден непарний тригонометричний поліном $T_{n-1}(t)$ порядку $n-1$ не може інтерполювати будь-яку функцію 
\begin{equation} \label{GrindEQ__39_} 
	K_m(\overline{\alpha},t)=\sum\limits^m_{i=1}{\alpha_i}D_{r_i}(t)=\sum\limits^m_{i=1}{\alpha_i(-1)^{\frac{r_i-1}{2}}}\sum\limits^\infty_{k=1}{\frac{{\sin k}t}{k^{r_i}}}, r_i=2l_i-1, l_i\in {\mathbb N},                
\end{equation} 
на інтервалі $\left(0,2\pi\right)$ більше ніж в $2(n+m)-3$ точках. 

    Не зменшуючи загальності, можемо вважати, що непарні різні числа $r_i$ в наборі $\overline{r}=(r_1,r_2,..,r_m)$ в \eqref{GrindEQ__38_} і \eqref{GrindEQ__39_} є розташованими у зростаючому порядку, тобто є такими, що $1\le r_1<r_2<...<r_m$.

   Далі нам знадобиться деяка додаткова інформація про ядра $D_r(t)$ вигляду \eqref{GrindEQ__5_}.
    Ядро $D_r(t),r\in {\mathbb N}{\rm ,}$ на періоді $(0,2\pi)$  при $r\ge 2$ є неперервно диференційовним до порядку $r-2$ включно $\left(D_r\in C^{r-2},r\ge 2\right)$ і істотно обмеженим при $r=1 (D_1\in L_{\infty })$. Крім того, для похідних $D^{(j)}_r$ порядку $j$ виконуються рівності
\begin{equation} \label{GrindEQ__40_} 
	D^{(j)}_r(t)=D_{r-j}(t), j=\overline{1,r-2}, t\in \left[0,2\pi\right], 
\end{equation} 
\begin{equation} \label{GrindEQ__41_} 
	D^{(r-1)}_r(t)=D_1(t)=\frac{\pi-t}{2}, t\in \left(0,2\pi\right). 
\end{equation} 
(див., наприклад, [11, c. 59-60]).

    На інтервалі $\left(0,2\pi\right)$ ядра $D_r(t)$ є алгебраїчними многочленами порядку $r$  і пов'язані з класичними многочленами Бернуллі $B_r(t)$  за допомогою рівності  
\begin{equation} \label{GrindEQ__42_} 
	D_r(t)=-\frac{(2\pi)^n}{2\cdot n!}B_r\left(\frac{t}{2n}\right), t\in \left(0,2\pi\right),r\in {\mathbb N}, 
\end{equation} 
де     $B_n(x)=\sum\limits^n_{k=0}{C^k_nB_{n-k}x^k},x\in \left(0,1\right),n\in {\mathbb N}{\rm ,}$  $C^k_n=\left( \begin{array}{c}
	n \\ 
	k \end{array}
\right)=\frac{n!}{k!(n-k)!}$ -- біноміальні коефіцієнти, $B_n$ -- числа Бернуллі $\left(B_n=-\frac{1}{n+1}\sum\limits^n_{k=1}{C^{k+1}_{n+1}B_{n-k},n\in {\mathbb N}, B_0=1}\right)$.

    Із наведених властивостей ядер Бернуллі випливає, що для будь-якої нетривіальної лінійної комбінації  $K_m({\overline{\alpha}}^*,t)$ вигляду \eqref{GrindEQ__39_} і довільного непарного тригонометричного полінома $T_{n-1}$ порядку не вищого  за $n-1$  різниця 
\begin{equation} \label{GrindEQ__43_} 
	f(t)=K_m\left(\overline{\alpha},t\right)-T_{n-1}(t) 
\end{equation} 
може мати на $[0,2\pi)$ лише скінченне число ізольованих нулів ${\left\{t_i\right\}}^N_{i=0},$
\[t_0=0<t_1<t_2<...<t_N=2\pi.\] 
Для кожного ізольованого нуля $t_i$ функції $f\in C$ існує достатньо мале $\varepsilon>0$ таке, що $f(t_i-\tau)f(t_i+\tau)\neq0,$ $\forall \tau\in (0,\varepsilon)$. 
 При цьому, якщо ${\rm sign}\left(f(t_i-\tau)f(t_i+\tau)\right)=-1$, то такий нуль називають простим нулем функції $f$, а якщо ${\rm sign}\left(f(t_i-\tau)f(t_i+\tau)\right)=1$, то такий нуль називають подвійним нулем функції $f$.

Нехай $Z(f)=Z\left(f;[0,2\pi)\right)$ підраховує кількість нулів неперервної функції на періоді $[0,2\pi)$ таким чином, що прості нулі враховуються один раз, а подвійні нулі -- двічі, тобто $Z\left(f;[0,2\pi)\right)=p+2d,$ де $p$ -- число простих нулів, а $d-$число подвійних нулів на періоді $[0,2\pi)$. Аналогічно означається і функція числа нулів  на інтервалі $\left(0,2\pi\right)$ (а також на будь-якому з відрізків $\left[a,b\right],\left(a,b\right),$
 $[\left.a,b\right),(a,b]$, включених в $\left(0,2\pi\right)$, в припущенні, що $f$ є неперервною на $\left(0,2\pi\right)$).

    Наша подальша мета полягає в тому, щоб показати, що для функції $f(t)$ вигляду \eqref{GrindEQ__39_} має місце нерівність
\begin{equation} \label{GrindEQ__44_} 
	Z\left(f;\left(0,2\pi\right)\right)\le 2n+2m-3.                                                  
\end{equation} 
Припустимо, що $Z\left(f;\left(0,2\pi\right)\right)=s$. Оскільки $f=f^{(0)}$ непарна функція, то $Z\left(f;(0,2\pi)\right)=2s+1$. Якщо $r_1>1$ , то для всіх похідних $f^{(j)},j=\overline{0,r_1-2}$, в силу \eqref{GrindEQ__40_} виконується включення $f^{(j)}\in C$ і в результаті послідовного застосування теореми Ролля та врахування того, що неперервна $2\pi$-періодична функція має парне число нулів, одержуємо             
\begin{equation} \label{GrindEQ__45_} 
	Z\left(f^{(j)};[\left.0,2\pi\right)\right)\ge 2s+2, j=0,1,...,r_1-2. 
\end{equation} 
При $j=r_1-1$ в силу \eqref{GrindEQ__41_} функція $f^{(j)}=f^{(r_1-1)}$ має розрив в точці 0 і в силу теореми Ролля і \eqref{GrindEQ__45_}  
\[Z\left(f^{(r_1-1)};(0,2\pi)\right)\ge 2s+1.\] 
Тоді для парної функції $f^{(r_1)}$ в силу \eqref{GrindEQ__39_} і \eqref{GrindEQ__41_} на $\left(0,2\pi\right)$ 
\begin{equation} \label{GrindEQ__46_} 
	f^{(r_1)}(t)=\sum\limits^m_{k=2}{\alpha_iD_{r_i-r_1}}(t)-T^{(r_1)}_{n-1}(t)-\frac{\alpha_1}{2},                                         \end{equation} 
і оскільки права частина в \eqref{GrindEQ__46_} є неперервною на всьому періоді $[0,2\pi)$ функцією, то в силу теореми Ролля для доозначеної за неперервністю в точці $0$ функції $f^{(r_1)}$ маємо  
\begin{equation} \label{GrindEQ__47_} 
	Z\left(f^{(r_1)};\left(0,2\pi\right)\right)\ge2s.                                                         
\end{equation} 
Аналогічно для всіх похідних $f^{(j)},j=r_1+1,r_1+2,...,r_2-2$, в силу \eqref{GrindEQ__40_} виконуються включення $f^{(j)}\in C$ і в результаті послідовного застосування теореми Ролля та врахування того, що всяка неперервна на періоді функція має парне число нулів, одержуємо
\begin{equation} \label{GrindEQ__48_} 
	Z\left(f^{(j)};[\left.0,2\pi\right)\right)\ge2s, j=r_1+1,r_1+2,...,r_2-2. 
\end{equation} 
При $j=r_2-1$ в силу \eqref{GrindEQ__41_} функція $f^{(j)}=f^{(r_2-1)}$ має розрив в точці $0$ і в силу теореми Ролля  
\[Z\left(f^{(r_2-1)};(0,2\pi)\right)\ge2s-1.                                             \] 
Тоді для парної функції $f^{(r_2)}$ в силу \eqref{GrindEQ__39_} і \eqref{GrindEQ__41_} на $\left(0,2\pi\right)$ 
\begin{equation} \label{GrindEQ__49_} 
	f^{(r_2)}(t)=\sum\limits^m_{k=3}{\alpha_iD_{r_i-r_2}}(t)-T^{(r_2)}_{n-1}(t)-\frac{\alpha_2}{2},                                    
\end{equation} 
і оскільки права частина в \eqref{GrindEQ__49_} є неперервною на всьому періоді $[\left.0,2\pi\right)$, то в силу теореми Ролля для доозначеної за неперервністю в точці $0$ функції $f^{(r_2)}$ маємо 
\begin{equation} \label{GrindEQ__50_} 
	Z\left(f^{(r_2)};[\left.0,2\pi\right)\right)\ge2(s-1).                                           
\end{equation} 
Аналогічно для всіх похідних $f^{(j)},j=r_2+1,r_2+2,...,r_3-2$, в силу \eqref{GrindEQ__40_} виконуються включення $f^{(j)}\in C$ і в результаті послідовного застосування теореми Ролля та  парності числа нулів  неперервної періодичної функції,  одержуємо
\begin{equation} \label{GrindEQ__51_} 
	Z\left(f^{(j)};[\left.0,2\pi\right)\right)\ge2(s-1), j=r_2+1,r_2+2,...,r_3-2. 
\end{equation} 
Продовжуючи міркування тим же чином, одержимо, що при $j=r_{m-1}$ для парної функції $f^{(r_{m-1})}$ в силу \eqref{GrindEQ__39_} і \eqref{GrindEQ__41_} на інтервалі $\left(0,2\pi\right)$
\begin{equation} \label{GrindEQ__52_} f^{(r_{m-1})}(t)=\alpha_mD_{r_m-r_{m-1}}(t)-T^{(m-1)}_{n-1}(t)-\frac{\alpha_{m-1}}{2}
\end{equation} 
     і оскільки права частина в \eqref{GrindEQ__52_} є неперервною на періоді $[0,2\pi)$, то для доозначеної за неперервністю в тоці 0 функції $f^{(r_{m-1})}$ в силу теореми Ролля
\begin{equation} \label{GrindEQ__53_} 
	Z\left(f^{(r_{m-1})};[0,2\pi)\right)\ge2(s-m+2). 
\end{equation} 

Тоді для всіх похідних $f^{(j)},j=r_{m-1}+1,r_{m-1}+2,...,r_{m-1}-2$, в силу \eqref{GrindEQ__40_} виконуються включення $f^{(j)}\in C$ і в результаті послідовного застосування \eqref{GrindEQ__33_}, теореми Ролля та врахування парності числа нулів функцій із $C$ матимемо 
\begin{equation} \label{GrindEQ__54_} 
	Z\left(f^{(j)};[0,2\pi)\right)\ge2(s-m+2), j=r_{m-1}+1,r_{m-1}+2,...,r_{m}-2. 
\end{equation} 
При $j=r_m-1$ функція $f^{(j)}=f^{(r_m-1)}$  в силу \eqref{GrindEQ__41_} має розрив в точці 0 і в силу \eqref{GrindEQ__54_} і теореми Ролля            
\begin{equation} \label{GrindEQ__55_} 
	Z\left(f^{(r_m-1)};(0,2\pi)\right)\ge2s-2m+3 
\end{equation} 
Парна функція $f^{(r_m)}$ в силу \eqref{GrindEQ__39_} і \eqref{GrindEQ__41_} на $\left(0,2\pi\right)$ має вигляд $f^{(r_m)}(t)=-T^{(m)}_{n-1}(t)-\frac{\alpha_m}{2}={\overline{T}}_{n-1}(t)$, тобто є тригонометричним поліномом  ${\overline{T}}_{n-1}$ порядку не вищого ніж $n-1$, а, отже, на $\left(0,2\pi\right)$ має не більше за $2n-2$ коренів, тобто 
\begin{equation} \label{GrindEQ__56_} 
	Z\left(f^{(r_m)};(0,2\pi)\right)\le2n-2. 
\end{equation} 
З іншого боку, в силу теореми Ролля із \eqref{GrindEQ__55_} одержуємо
\begin{equation} \label{GrindEQ__57_} 
	Z\left(f^{(r_m)};(0,2\pi)\right)\ge2s-2n+2. 
\end{equation} 
Припустимо, що нерівність \eqref{GrindEQ__44_} місця не має, тобто, що $Z\left(f;(0,2\pi)\right)=2s+1>2n+2m-3$, а, отже, $2s>2(n+m-2)$. Тоді в силу \eqref{GrindEQ__57_} отримуємо нерівність $Z\left(f^{(r_m)};(0,2\pi)\right)>2n-2$, яка суперечить нерівності \eqref{GrindEQ__56_}. Одержана суперечність доводить істинність формули \eqref{GrindEQ__44_}.

   Із нерівності \eqref{GrindEQ__44_}, застосованої до функції
\begin{equation} \label{GrindEQ__58_} 
	f^*(t)=K_m({\overline{\alpha}}^*,t)-T^*_{n-1}(t), 
\end{equation} 
де $K_m({\overline{\alpha}}^*,t)$  має вигляд \eqref{GrindEQ__30_}, та $T^*_{n-1}$  такий, що виконуються рівності \eqref{GrindEQ__37_}, при \mbox{$k=\overline{1,2(n+m)-3}$}  випливає, що всі нулі $t_k$ вигляду \eqref{GrindEQ__36_} для різниці $f^*(t)$ вигляду \eqref{GrindEQ__58_} є простими нулями, в яких вона почережно змінює знак. Теорему 3 доведено.

     \textbf{Доведення теореми 4}. Оскільки спряжені ядра Пуассона   
\begin{equation} \label{GrindEQ__59_} 
	P_{q_i,1}(t)=\sum\limits^\infty_{k=1}{q^k_i{\sin k}t=\frac{q_i{\sin t}}{1-2q_i{\cos t\ }+q^2_i}}, q_i\in (0,1), i=\overline{1,m}, 
\end{equation} 
є неперервними на $[0,2\pi)$ функціями, то відповідно до наслідку 1 існує нетривіальна лінійна комбінація $K_m({\overline{\alpha}}^*,t)$ вигляду \eqref{GrindEQ__32_} і непарний тригонометричний поліном $T^*_{n-1}(t)$ порядку не вищого за $n-1$, який інтерполює $K_m({\overline{\alpha}}^*,t)$ на $\left(0,2\pi\right)$ в $2(n+m)-3$ точках вигляду \eqref{GrindEQ__36_}. 

    Крім того, в силу непарності і неперервності на $\left[0,2\pi\right]$ різниця $K_m\left({\overline{\alpha}}^*,t\right)-T^*_{n-1}(t)$ дорівнює 0 і при $t=0$. Тобто рівність \eqref{GrindEQ__37_} виконується на періоді $[0,2\pi)$ при всіх $k=\overline{0,2(n+m)-3}$, а, отже, число точок інтерполяції по меншій мірі $2(n+m-1)$ . 

    З іншого боку неважко показати, що жоден непарний тригонометричний поліном $T_{n-1}(t)$ порядку $n-1$ не може інтерполювати будь-яку функцію 
\begin{equation} \label{GrindEQ__60_} 
	K_m(\overline{\alpha},t)=\sum\limits^m_{i=1}{\alpha_iP_{q_i,1}}(t), \alpha_i\in {\mathbb R}, q_i\in (0,1), 
\end{equation} 
на півінтервалі $[0,2\pi)$ більше, ніж в $2(n+m-1)$ точках. Дійсно, в силу \eqref{GrindEQ__59_} і \eqref{GrindEQ__60_}
\[K_m(\overline{\alpha},t)=P^{\overline{\alpha},\overline{q}}_{m,1}(t)=\sum\limits^m_{i=1}{\alpha_i}\frac{q_i{\sin t }}{1-2q_i{\cos t}+q^2_i}=\] 
\begin{equation} \label{GrindEQ__61_} 
	=\frac{\sin t}{\prod\limits_{i=1}^{m}\left(1-2q_{i} \cos t+q_{i}^{2} \right) }  \sum _{i=1}^{m}\alpha _{i} q_{i}  \prod _{\begin{array}{c} {k=1} \\ {k\ne i} \end{array}}^{m}\left(1-2q_{k} \cos t+q_{k}^{2} \right) =\frac{{\overline{T}}_m(t)}{\prod\limits^m_{i=1}{\left(1-q_i{\cos t }+q^2_i\right)}}, 
\end{equation} 
де ${\overline{T}}_m(t)-$ деякий непарний тригонометричний многочлен порядку не вищого за $m$. В силу \eqref{GrindEQ__61_}  рівність  $K_m(\overline{\alpha},t)-T_{n-1}(t)=0$ можлива лише при виконанні співвідношення
\begin{equation} \label{GrindEQ__62_} 
	{\overline{T}}_m(t)-T_{n-1}(t)\prod^m_{i=1}{\left(1-2q_i{\cos t }+q^2_i\right)}=0.                                      
\end{equation} 
Ліва частина рівності \eqref{GrindEQ__62_} є тригонометричним многочленом порядку не вище $m+n-1$, а тому може мати на періоді не більше за $2(n+m-1)$ нулів, з урахуванням їх кратності.

    Із сказаного випливає, що всі нулі $t_k=\frac{k?
    \pi}{n+m-1}, k=\overline{0,2(n+m)-3}$ , для різниці $K_m({\overline{\alpha}}^*,t)-T^*_{n-1}(t)$ є простими нулями, в яких ця різниця почережно змінює знак. Теорему 4 доведено.

 \textbf{Конфлікт інтересів.} Автори заявляють, що вони не мають потенційного конфлікту інтересів щодо дослідження у цій статті.

 \textbf{Фінансування}. Цю роботу частково підтримано грантом H2020-MSCA-RISE-2019, проект № 873071 (SOMPATY: spectral optimization: from mathematics to physics and advanced technology).  

 \textbf{Авторські внески}. Усі автори зробили рівний внесок у роботу.

                                               Література

\begin{enumerate}
	\item   Степанец А.И., \textit{Методы теории приближений}: В 2 ч. -- Киев: Ин-т математики НАН Украини, 2002. -- Ч. 1. -- 427 с.
	
	\item  Крейн М.Г., \textit{К теории наилучшего приближения периодических функций }Докл. АН СССР. -1938, т.18, № 4 -5. -- C.245-249.
	
	\item  Ахиезер Н.И., \textit{О наилучшем приближении одного класса непрерывных} \textit{периодических функций}, Докл. АН СССР. -- 1937, т.17, № 9. -- с.451-454.
	
	\item  Favard J., \textit{Sur l'approximation des fonctions p}$\acute{e}$\textit{riodiques par des polynomes} \textit{trigonom}$\acute{e}$\textit{triques}, C.r. Acad. Sci. -- 1936. -- 203. P. 1122-1124.
	
	\item Favard J.,  \textit{Sur les meilleurs proc}$\acute{e}$\textit{des d'approximations des certain classes des fonctions par des polynomes trigonom}$\acute{e}$\textit{triques}, Bull. de  Scienses  Math. -- 1937. -- 15. --P. 209-224, 243-256.
	
	\item  Ахиезер Н.И., Крейн М.Г., \textit{О наилучшем приближении тригонометрически-ми суммами дифференцируемых периодических функций}, Докл. АН СССР. -- 1937. -- 15. № 3. -- С. 107-112.
	
	\item  Никольский С.М., \textit{Приближение функций тригонометрическими} \textit{полиномами в среднем}, Изв. АН СССР. Сер. мат. -- 1946. -- 10. -- С. 207-256.
	
	\item  Nagy B., $\ddot{U}$\textit{ber gewisse Extremalfragen bei transformierten trigonometrischen} \textit{Entwicklungen 1. Periodischer Fall}, Berichte der math. --phys. Kl. Acad. Der Wiss. zu Leipzig. -- 1938. -- 90. --P. 103-134.
	
	\item  Бушанский А.В., \textit{О наилучшем в среднем гармоническом приближении} \textit{некоторых функций}, Исследования по теории приближения функций и их приложения. -- Киев: Ин-т математики АН УССР. -- 1978. -- С. 29-37.
	
	\item  Шевалдин В.Т., \textit{Поперечники классов сверток с ядром Пуассона}, Матем. заметки. -- 1993. -- 52, № 2. -- С. 145-151.
	
	\item  Корнейчук Н.П., \textit{Экстремальные задачи теории приближения} -- М.: Наука, Главная редакция физико-математической литературы, 1976. -- 320 с.
	
	\item  Pincus A., \textit{n-Widths in Approximation Theory}, Springer-Verlag Berlin Heidelberg, 1985. -- 294~p.
	
	\item  Корнейчук Н.П., \textit{Точные константы в теории приближения} -- М.: Наука, Главная редакция физико-математической литературы, 1987. -- 424 с.
	
	\item  Сердюк А.С., \textit{Поперечники та найкращі наближення класів згорток} \textit{періодичних функцій} -- Укр. мат. журн. -- 1999. -- Т. 51, № 5. -- С. 674-687.
	
	\item  Serdyuk A.S., Bodenchuk V.V., \textit{Exact values of Kolmogorov widths of classes of} \textit{Poisson integrals}, Jornal of Approximation Theory. -- 2013. -- 173, № 9. -- Р. 89-109.
	
	\item  Боденчук В.В., Сердюк А.С., \textit{Точні оцінки колмогоровських поперечників} \textit{класів аналітичних функцій І} , Укр. мат. журн. -- 2015. -- Т. 67, № 6. -- С. 719-738.
	
	\item  Боденчук В.В., Сердюк А.С., \textit{Точні оцінки колмогоровських поперечників} \textit{класів аналітичних функцій ІІ}, Укр. мат. журн. -- 2015. -- Т. 67, № 8. -- С. 1011-1018.
	
	\item  Grabova U.Z., Kal'chuk I.V., \textit{Approximation of the classes }$W^r_{\beta,\infty}$ \textit{ by three}-\textit{garmonic Poisson integrals}, Carpathion Math. Publ. -- 2019. -- T. 11, № 2. - P. 321-334.

	\item  Абдулаєв Ф.Г., Харкевич Ю.І., \textit{Наближення класів }$C^\psi_\beta H^\omega$ \textit{ бігармонічними} \textit{інтегралами Пуассона}, Укр. мат. журн. -- 2020. -- Т. 72, № 1. -- С. 20-35.
	
	\item  Сердюк А.С., \textit{Найкращі наближення і поперечники класів згорток періодичних функцій високої гладкості}, Укр. мат. журн. -- 2005. -- Т. 57, № 7. -- С. 946-971.
	
	\item  Степанец А.И., \textit{Методы теории приближений:} в 2 ч. -- Киев: Ин-т математики НАН Украины, 2002. -- Ч. ІІ. -- 468 с.
	
	\item  Сердюк А.С., \textit{Про найкраще наближення на класах згорток періодичних} \textit{функцій}, Теорія наближення функцій та суміжні питання/ Праці Інституту-ту математики НАН України. -- 35. -- К.: Ін-т математики НАН України. 2002. -- С. 172-194.
\end{enumerate}

\end{document}